\documentclass[12pt,a4paper]{article}
\usepackage{amsfonts,amssymb,amsmath,amsthm,cases}
\usepackage[colorlinks=true,linkcolor=blue,citecolor=red, breaklinks = true]{hyperref}
\newtheoremstyle{thmstyle}{}{}{}{}{\bf}{}{ }{}
\theoremstyle{thmstyle}
\numberwithin{equation}{section}

\newtheorem{thm}{Theorem}[section]

\allowdisplaybreaks
\begin{document}
\begin{center}
{\bf {On Generalized Free Electron Laser Equation}}\\
Dharmendra Kumar Singh\\
Department of Mathematics\\
School of Basic Sciences\\
CSJM University, Kanpur\\
(UP)India\\
drdksinghabp@gmail.com\\
\end{center}
\begin{abstract}
The aim of this paper is to derive a solution for a generalized free electron laser equation in terms of the incomplete Mittag-Leffler function and in terms of the incomplete Wright function..
\end{abstract}
\section{Introduction}
Fractional Calculus is one of the most important tools for theoretical and experimental physics. Fractional calculus is of great scientific significance as it bestows a unified approach to various physical problems and yields fundamental exploratory ideas.In view of the effect of fractional calculus, it has been used by many researchers \cite{2,3,4,6,7,8, 9,10,11,12,13,15,16} in the study of free electron lasers. Dattoli et al. \cite{6, 9} described the unsaturated behavior of the free electron laser in the first-order integro-differential equation of Volterra type as
\begin{equation}\label{1.1}
D_{\mu}h(\mu)=-i\pi g_{0}\int_{0}^{\mu}\psi h(\mu-\psi)e^{i\nu\psi}d\psi,~~~~~ 0\leq\mu\leq 1.
\end{equation}
where $\mu$ is a dimensionless time variable, $g_{0}$ is a positive constant called the small-signal gain and the detuning parameter is the constant $\nu$. The function $h(\mu)$ is the complex-field amplitude which is assumed to be dimensionless satisfying the initial condition $h(0)=1$. The exact closed-form solution of (\ref{1.1}) valid in the region of practical interest and suitable for numerical computation is given in \cite{8}.\\
The fractional integral of an arbitrary order is the generalization of the n-th order integral$(n\in N)$. The most fundamental definition of the fractional integral of arbitrary order $a,$ due to Riemann and Liouville \cite{1}, is given by
\begin{equation}\label{1.2}
(R^{a}_{0, x}f)(x)=(I^{a}_{0,x}f)(x)=(D^{-a}_{0, x}f)(x)=\frac{1}{\Gamma(a)}\int_{0}^{x} (x-t)^{a-1} f(t) dt
\end{equation}
\begin{equation}\label{1.3}
=\frac{x^{a}}{\Gamma(a)}\int_{0}^{1} (1-u)^{a-1} f(xu) du,~~~~(\Re(a)>0).
\end{equation}
and 
\begin{equation}\label{1.4}
(D^{a}_{0,x}f)(x)=D^{n}_{0,x}(D^{a-n}_{0,x}f)(x)=D^{n}_{0,x}(I^{n-a}_{0,x}f)(x).
\end{equation}
Laplace transform of the fractional derivative
\begin{equation}\label{1.5}
L\{D^{a}_{t};s\}=s^{a}F(s)-\sum_{r=1}^{n}s^{r-1}D^{a-r}_{t}f(t)\vert_{t}=0, ~~~~(r-1<a\leq r;r\in N).
\end{equation}
The Caputo derivatives of the casual function $f(t)$ (that is $f(t)=0$ for $t<0 $) with $a>0$ was defined by Caputo \cite{14} in the form
\begin{equation}\label{1.6}
_{0}^{c}D^{a}_{t}f(t)=(_{0}^{c}D^{a}_{t}f)(t)=\frac{d^{a}}{dx^{a}}f(t)=\frac{1}{\Gamma(m-a)}\int_{0}^{\tau}\frac{f^{(m)}(u)}{(t-u)^{a-m+1}}du,
\end{equation}
where $m-1<a<m, m\in N$ and $f^{(m)}(t)$ denote the usual derivative of $f(t)$ of order $m\in N_{0}$.\\
The Laplace transform of the Caputo derivative  in the form
\begin{equation}\label{1.7}
 L\{\frac{d^{a}}{dt^{a}}f(t); s\}=s^{a} F(s)-\sum_{r=0}^{m-1}s^{a-r-1}f^{(r)}(0); m-1<a\leq m; m\in N,
\end{equation}
where $F(s)$ is the Laplace transform of $f(t)$
\begin{equation}\label{1.8}
L\{f(t)\}=F(s)=\int_{0}^{\infty}e^{-st}f(t)dt,
\end{equation}
where $\Re(s)>0$, which may be symbolically written as
\begin{equation*}
F(s)=L\{f(t); s\}\text{(or)} f(t)=L^{-1}\{F(s); t\}.
\end{equation*}
The homogenous form of equation (\ref{1.1}) is given by Boyadjiev et al. \cite{12}, which is as follows
\begin{equation}\label{1.9}
D^{a}_{\mu}h(\mu)=\lambda\int_{0}^{\mu}\psi h(\mu-\psi)e^{i\nu\psi}d\psi+be^{i\mu},~~~~~ 0\leq\mu\leq 1,
\end{equation}
where $b, \lambda\in C$ and $\nu\in R.$\\
The generalization form of the equation (\ref{1.2}) was presented by Al-Shammery et al.\cite{3} in the following manner
\begin{equation}\label{1.10}
D^{a}_{\mu}h(\mu)=\lambda\int_{0}^{\mu}\psi^{\delta} h(\mu-\psi)e^{i\nu\psi}d\psi+be^{i\nu\mu},~~~~~ 0\leq\mu\leq 1,
\end{equation}
where $b, \lambda\in C$ $\nu\in R.$and $\delta>-1$\\
Al-Shammery et al. \cite{4} established a new form of the equation (\ref{1.3}) by including function $\Phi(.)$
\begin{equation}\label{1.11}
D^{a}_{\mu}h(\mu)=\lambda\int_{0}^{\mu}\psi^{\delta} h(\mu-\psi)\Phi(\beta, \delta+1;i\nu\psi) d\psi+b\Psi(\beta', 1; i\nu\psi),~~~~~ 0\leq\mu\leq 1,
\end{equation}
where $b, \lambda\in C$ $\nu, \beta, \beta'\in R, a>0$ and $\delta>-1$.
Now, in this paper we will deal with equation (\ref{1.11}) with the function introduced by Singh and Porwal \cite{5} known as the incomplete Mittag-Leffler function and whose definition is as follows\\
 \begin{equation}\label{1.12}
 E_{a, b}^{[\delta, x]}=\sum_{k=0}^{\infty}\frac{[\delta, x]_{k}}{\Gamma(a k+b)}\frac{z^{k}}{k!}
 \end{equation}
 \begin{equation}\label{1.13}
  E_{a, b}^{(\delta, x)}=\sum_{k=0}^{\infty}\frac{(\delta, x)_{k}}{\Gamma(a k+b)}\frac{z^{k}}{k!},
\end{equation}
  where $a, b, \delta\in C$; $\Re(a)>0, \Re(b)>0, \Re(\delta)>0$ and $[\delta, x]_{k}$ and $(\delta, x)_{k}$ represent incomplete Pocchammer symbol and these incomplete Pocchammer symbols satisfy the following decomposition
 \begin{equation}\label{1.14}
 (\lambda; x)_{\nu}+[\lambda; x]_{\nu}=(\lambda)_{\nu}~~~~~~~~(\lambda, \nu\in C; x\geq 0).
 \end{equation}
  {\bf{Incomplete Wright Function:}}\\
 Incomplete generalized hypergeometric function \cite{5} is defined as
\begin{equation*}
{_{q}{\overline{\Psi}}_{q} (z)}={_{q}{\overline{\Psi}}_{q}}\left[\begin{array}{cc} [a_{1}, \alpha_{1}, x]...(a_{p}, \alpha_{p})\\(b_{1}, \beta_{1}, x]...(b_{q}, \beta_{q})\end{array}\vert z\right]
\end{equation*}
\begin{equation}\label{1.15}
=\sum_{k=0}^{\infty}\frac{\Gamma(a_{1}+\alpha_{1}k, x) \Gamma(a_{2}+\alpha_{2}k)...\Gamma(a_{p}+\alpha_{p}k)}
{\Gamma(b_{1}+\beta_{1}k, x) \Gamma(b_{2}+\beta_{2}k)...\Gamma(b_{q}+\beta_{q}k)}\frac{z^{k}}{k!}
\end{equation}
\begin{equation*}
{_{q}{\underline{\Psi}}_{q} (z)}={_{q}{\underline{\Psi}}_{q}}\left[\begin{array}{cc} (a_{1}, \alpha_{1}, x)...(a_{p}, \alpha_{p})\\(b_{1}, \beta_{1}, x]...(b_{q}, \beta_{q})\end{array}\vert z\right]
\end{equation*}
 \begin{equation}\label{1.16}
 \sum_{k=0}^{\infty}\frac{\gamma(a_{1}+\alpha_{1}k, x) \Gamma(a_{2}+\alpha_{2}k)...\Gamma(a_{p}+\alpha_{p}k)}
 {\gamma(b_{1}+\beta_{1}k, x) \Gamma(b_{2}+\beta_{2}k)...\Gamma(b_{q}+\beta_{q}k)}\frac{z^{k}}{k!}
 \end{equation}
\smallskip

With decomposition formula (\ref{1.14}) incomplete  Wright function (\ref{1.12}) and (\ref{1.13}) as classical Wright function $_{p}\Psi_{q}(z)$ \cite{1} and this particular function is an entire function if there hold the condition
\smallskip

  \begin{equation}\label{1.17}
\sum_{j=1}^{q}\beta_{j}-\sum_{i=1}^{p}\alpha_{j}>1.
  \end{equation}

\section{Theorems and Results}
\begin{thm} 
Let $b, c, \omega, \in C, \rho\in \Re,a>0, \rho>0,\Re(b)>0$ and $g(t)$  is assumed to be continuous on every finite interval 
$\left[0, T \right]$,$0<T<\infty$, and of the exponential order $e^{\rho\mu}$, when $\mu\to\infty$. Then for the generalized Free Electron Laser equation for the fractional integro-differential equation of Volterra  type
\begin{equation}\label{2.1}
D^{a}_{\mu} h(\mu)=\omega\int_{0}^{\mu} t^{b-1}  E^{[c; x]}_{\rho,  b}(i\zeta t^{\rho}) h(\mu-t) dt+\delta g(t)~~~~(0\leq\mu\leq 1).
\end{equation}
together with the initial conditions
\begin{equation}\label{2.2}
D^{a-r}_{\mu}h(\mu)\vert_{\mu}=0=b_{r}, r=1, ..., n=-[-\Re(a)], \left( n-1<a\leq n \right); n\in N,
\end{equation}
where $b_{1}, b_{2}, b_{3}, ...., b_{r}\in \Re$, there exist a unique continuous solution given by
\begin{equation}\label{2.3}
h(\mu)=\sum_{r=1}^{n}b_{r}y_{r}(\mu)+\delta\int_{0}^{\mu}\aleph(\mu-t)g(t)dt
\end{equation}
where
\begin{equation}\label{2.4}
y_{r}(\mu)=\sum_{k=0}^{\infty}\omega^{k}\mu^{a+(a+b)k-r}E^{[c; x]k}_{\rho,1+a+(a+b)k-r}(i\zeta \mu^{\rho})\\
\end{equation}
\begin{equation}\label{2.5}
y_{r}(\mu)=\sum_{k=0}^{\infty}\omega^{k}\mu^{a+(a+b)k-r}\frac{1}{\Gamma(ck)}{_{1}}\overline{\Psi}_{1}\left[\begin{array}{cc}[ck, 1; xk]\\(1+a+(a+b)k-r, \rho)\end{array}\vert i\zeta\mu^{\rho}\right]
\end{equation}
\begin{equation}\label{2.6}
\aleph(x)=\sum_{k=0}^{\infty}\omega^{k} x^{a+(a+b)k-1}E^{[c;x]k}_{\rho, a+(a+b)k}(i\zeta x^{\rho})
\end{equation}
\begin{equation}\label{2.7}
\aleph(x)=\sum_{k=0}^{\infty}\omega^{k} x^{a+(a+b)k-1}\frac{1}{\Gamma(ck)}{_{1}}\overline{\Psi}_{1}\left[\begin{array}{cc}[ck, 1; xk]\\(a+(a+b)k, \rho)\end{array}\vert i\zeta\mu^{\rho}\right]
\end{equation}
\end{thm}
\smallskip

\begin{proof}
Taking Laplace on both side of the integro-differential equation (\ref{2.1}) and using (\ref{1.12}) and (\ref{1.5}),  we find that
\begin{equation}\label{2.8}
L \{D^{a}_{\mu}h{\left(\mu\right)}\}=L\{\omega\int^{\mu}_{0} t^{b-1}h(\mu-t) E^{[c; x]}_{\rho,  b}(i\zeta t^{\rho})dt+\delta g(t)\}
\end{equation}
\begin{equation}\label{2.9}
s^{a}H(s)-\sum_{r=0}^{n}s^{r-1}D^{a-r}_{\mu}h(\mu)\vert_{\mu=0}=\omega s^{-b}(1-i\zeta s^{-\rho})^{-[c;x]} H(s)+\delta G(s)
\end{equation}
\begin{equation}\label{2.10}
\left[s^{a}-\omega s^{-b}(1-i\zeta s^{-\rho})^{-[c;x]} \right]H(s)=\sum_{r=0}^{n}s^{r-1}D^{a-r}_{\mu}h(\mu)\vert_{\mu=0}+\delta G(s)
\end{equation}
\begin{equation}\label{2.11}
\left[s^{a}-\omega s^{-b}(1-i\zeta s^{-\rho})^{-[c;x]} \right]H(s)=\sum_{r=0}^{n}s^{r-1}b_{r}+\delta G(s),
\end{equation}
where $D^{a-r}_{\mu}h(\mu)\vert_{\mu=0}$.
\begin{equation*}
H(s)=\sum_{r=1}^{n} s^{r-1} b_{r}\left[s^{a}-\omega s^{-b}(1-i\zeta s^{-\rho})^{-[c;x]}\right]^{-1}+\delta G(s)\left[s^{a}-\omega s^{-b}(1-i\zeta s^{-\rho})^{-[c;x]}\right]^{-1}
\end{equation*}
\begin{equation*}
H(s)=\sum_{r=1}^{n} s^{r-1} b_{r}s^{-a}\left[1-\omega s^{-a}s^{-b}(1-i\zeta s^{-\rho})^{-[c;x]}\right]^{-1}
\end{equation*}
\begin{equation*}
+\delta G(s)s^{-a}\left[1-\omega s^{-a}s^{-b}(1-i\zeta s^{-\rho})^{-[c;x]}\right]^{-1}
\end{equation*}
\begin{equation*}
=\sum_{r=1}^{n} s^{r-1} b_{r}s^{-a}\sum_{k=0}^{\infty}\left(\omega s^{-a}s^{-b}\right)^{k}(1-i\zeta s^{-\rho})^{-[c;x]k}
\end{equation*}
\begin{equation}\label{2.12}
+\delta G(s)s^{-a}\sum_{k=0}^{\infty}\left(\omega s^{-a}s^{-b}\right)^{k}(1-i\zeta s^{-\rho})^{-[c;x]k}
\end{equation}
\smallskip

where $\vert\omega s^{-a}s^{-b}(1-i\zeta s^{-\rho})^{-[c;x]}\vert<1$. Taking the inverse Laplace  transform on both sides of (\ref{2.12}), we get
\smallskip

\begin{equation*}
L^{-1}\{H(s)\}=L^{-1}\{\sum_{r=1}^{n} s^{r-1} b_{r}s^{-a}\sum_{k=0}^{\infty}\left(\omega s^{-a}s^{-b}\right)^{k}(1-i\zeta s^{-\rho})^{-[c;x]k}\}
\end{equation*}
\begin{equation*}
+L^{-1}\{\delta G(s)s^{-a}\sum_{k=0}^{\infty}\left(\omega s^{-a}s^{-b}\right)^{k}(1-i\zeta s^{-\rho})^{-[c;x]k}\}
\end{equation*}
\begin{equation*}
h(s)=L^{-1}\{\sum_{r=1}^{n} s^{r-1} b_{r}s^{-a}\sum_{k=0}^{\infty}\left(\omega s^{-a}s^{-b}\right)^{k}\sum_{n=0}^{\infty} (i\zeta s^{-\rho})^{n} ([c;x]k)_{n}\frac{1}{n!}\}
\end{equation*}
\begin{equation*}
+L^{-1}\{\delta G(s)s^{-a}\sum_{k=0}^{\infty}\left(\omega s^{-a}s^{-b}\right)^{k}\sum_{n=0}^{\infty}(i\zeta s^{-\rho})^{n}([c;x]k)_{n}\frac{1}{n!}\}
\end{equation*}
\begin{equation*}
h(s)=\sum_{r=1}^{n}b_{r}\sum_{k=0}^{\infty}\omega^{k}\sum_{n=0}^{\infty} (i\zeta)^{n} ([c;x]k)_{n}\frac{1}{n!}L^{-1}\{s^{r-a-(a+b)k-\rho n}\}
\end{equation*}
\begin{equation}\label{2.13}
+\delta \sum_{k=0}^{\infty}\omega^{k}\sum_{n=0}^{\infty}(i\zeta)^{n}([c;x]k)_{n}\frac{1}{n!}L^{-1}\{G(s)s^{-a-(a+b)k-\rho n}\}
\end{equation}
\smallskip

Using $L^{-1}\{\frac{1}{s^{n}}\}=\frac{t^{n-1}}{(n-1)!}$ in the above expression, then after some calculation, we arrive at
\begin{equation*}
=\sum_{r=1}^{n} b_{r}\sum_{k=0}^{\infty}\omega^{k}\mu^{a+(a+b)k-r}E^{[c; x]k}_{\rho,1+a+(a+b)k-r}(i\zeta \mu^{\rho})
\end{equation*}
\begin{equation}\label{2.14}
+\delta\int_{0}^{\mu}\sum_{k=0}^{\infty}\omega^{k} (\mu-t)^{a+(a+b)k-1}E^{[c;x]k}_{\rho, a+(a+b)k}(i\zeta (\mu-t)^{\rho})g(t)dt.
\end{equation}
\smallskip

Which is the required result and by using (\ref{1.15}) in (\ref{2.14}), we can achieve (\ref{2.5}) and (\ref{2.7}).\\
\end{proof}
\smallskip

\begin{thm}
Let $b, c, \omega, \in C, \rho\in \Re,a>0, \rho>0,\Re(b)>0$ and $g(t)$  is assumed to be continuous on every finite interval 
$\left[0, T \right]$,$0<T<\infty$, and of the exponential order $e^{\rho\mu}$, when $\mu\to\infty$. Then for the generalized Free Electron Laser equation for the fractional integro differential equation of Voltera type
\smallskip

\begin{equation}\label{2.15}
\frac{d^{a}}{d\mu^{a}} h(\mu)=\omega\int_{0}^{\mu} t^{b-1}  E^{[c; x]}_{\rho,  b}(i\zeta t^{\rho}) h(\mu-t) dt+\delta g(t)~~~~(0\leq\mu\leq 1).
\end{equation}
together with the initial conditions
\begin{equation}\label{2.16}
\frac{d^{a}}{d\mu^{a}} h(\mu)\vert_{\mu}=0=a_{r}, r=1, ..., n-1; \left( n-1<a\leq n \right); n\in N,
\end{equation}
where $a_{1}, a_{2}, a_{3}, ...., a_{n-1}\in \Re$, there exist a unique continuous solution given by
\begin{equation}\label{2.17}
h(\mu)=\sum_{r=0}^{n-1}a_{r}y_{r}(\mu)+\delta\int_{0}^{\mu}\aleph(\mu-t)g(t)dt
\end{equation}
where
\begin{equation}\label{2.18}
y_{r}(\mu)=\sum_{k=0}^{\infty}\omega^{k}\mu^{r+(a+b)k}E^{[c; x]k}_{\rho,1+(a+b)k+r}(i\zeta \mu^{\rho})
\end{equation}
\begin{equation}\label{2.19}
y_{r}(\mu)=\sum_{k=0}^{\infty}\omega^{k}\mu^{(a+b)k+r}\frac{1}{\Gamma(ck)}{_{1}}\overline{\Psi}_{1}\left[\begin{array}{cc}[ck, 1; xk]\\(1+(a+b)k+r, \rho)\end{array}\vert i\zeta\mu^{\rho}\right]
\end{equation}
\begin{equation}\label{2.20}
\aleph(x)=\sum_{k=0}^{\infty}\omega^{k} x^{a+(a+b)k-1}E^{[c;x]k}_{\rho, a+(a+b)k}(i\zeta x^{\rho})
\end{equation}
\begin{equation}\label{2.21}
\aleph(x)=\sum_{k=0}^{\infty}\omega^{k} x^{a+(a+b)k-1}\frac{1}{\Gamma(ck)}{_{1}}\overline{\Psi}_{1}\left[\begin{array}{cc}[ck, 1; xk]\\(a+(a+b)k, \rho)\end{array}\vert i\zeta\mu^{\rho}\right]
\end{equation}
\end{thm}
\smallskip

\begin{proof}
Taking Laplace on both side of the integro-differential equation (\ref{2.15}) and using (\ref{1.12}) and (\ref{1.5}), we find that
\begin{equation*}
L \{\frac{d^{a}}{d\mu^{a}}h{\left(\mu\right)}\}=L\{\omega\int^{\mu}_{0} t^{b-1}h(\mu-t) E^{[c; x]}_{\rho,  b}(i\zeta t^{\rho})dt+\delta g(t)\}.
\end{equation*}
\begin{equation*}
s^{a}H(s)-\sum_{r=0}^{n-1}s^{a-r-1}\frac{d^{a}}{d\mu^{a}} h(\mu)\vert_{\mu=0}=\omega s^{-b}(1-i\zeta s^{-\rho})^{-[c;x]} H(s)+\delta G(s).
\end{equation*}
\begin{equation*}
\left[s^{a}-\omega s^{-b}(1-i\zeta s^{-\rho})^{-[c;x]}\right]H(s)=\sum_{r=0}^{n-1}s^{a-r-1}\frac{d^{a}}{d\mu^{a}}h(\mu)\vert_{\mu=0}+\delta G(s).
\end{equation*}
\begin{equation*}
\left[s^{a}-\omega s^{-b}(1-i\zeta s^{-\rho})^{-[c;x]} \right]H(s)=\sum_{r=0}^{n-1}s^{a-r-1}b_{r}+\delta G(s),
\end{equation*}
where $\frac{d^{a}}{d\mu^{a}}h(\mu)\vert_{\mu=0}$.
\begin{equation*}
H(s)=\sum_{r=0}^{n-1} s^{a-r-1} b_{r}\left[s^{a}-\omega s^{-b}(1-i\zeta s^{-\rho})^{-[c;x]}\right]^{-1}+\delta G(s)\left[s^{a}-\omega s^{-b}(1-i\zeta s^{-\rho})^{-[c;x]}\right]^{-1}
\end{equation*}
\begin{equation*}
H(s)=\sum_{r=0}^{n-1} s^{-r-1} b_{r}\left[1-\omega s^{-a}s^{-b}(1-i\zeta s^{-\rho})^{-[c;x]}\right]^{-1}
\end{equation*}
\begin{equation*}
+\delta G(s)s^{-a}\left[1-\omega s^{-a}s^{-b}(1-i\zeta s^{-\rho})^{-[c;x]}\right]^{-1}
\end{equation*}
\begin{equation*}
=\sum_{r=0}^{n-1} s^{-r-1} b_{r}\sum_{k=0}^{\infty}\left(\omega s^{-a}s^{-b}\right)^{k}(1-i\zeta s^{-\rho})^{-[c;x]k}
\end{equation*}
\begin{equation}\label{2.22}
+\delta G(s)s^{-a}\sum_{k=0}^{\infty}\left(\omega s^{-a}s^{-b}\right)^{k}(1-i\zeta s^{-\rho})^{-[c;x]k}
\end{equation}
\smallskip

where $\vert\omega s^{-a}s^{-b}(1-i\zeta s^{-\rho})^{-[c;x]}\vert<1$. Taking the inverse Laplace transform on both sides of (\ref{2.22}), we get
\begin{equation*}
L^{-1}\{H(s)\}=L^{-1}\{\sum_{r=0}^{n-1} s^{-r-1} b_{r}\sum_{k=0}^{\infty}\left(\omega s^{-a}s^{-b}\right)^{k}(1-i\zeta s^{-\rho})^{-[c;x]k}\}
\end{equation*}
\begin{equation*}
+L^{-1}\{\delta G(s)s^{-a}\sum_{k=0}^{\infty}\left(\omega s^{-a}s^{-b}\right)^{k}(1-i\zeta s^{-\rho})^{-[c;x]k}\}
\end{equation*}
\begin{equation*}
h(s)=L^{-1}\{\sum_{r=0}^{n-1} s^{-r-1} b_{r}\sum_{k=0}^{\infty}\left(\omega s^{-a}s^{-b}\right)^{k}\sum_{n=0}^{\infty} (i\zeta s^{-\rho})^{n} ([c;x]k)_{n}\frac{1}{n!}\}
\end{equation*}
\begin{equation*}
+L^{-1}\{\delta G(s)s^{-a}\sum_{k=0}^{\infty}\left(\omega s^{-a}s^{-b}\right)^{k}\sum_{n=0}^{\infty}(i\zeta s^{-\rho})^{n}([c;x]k)_{n}\frac{1}{n!}\}
\end{equation*}
\begin{equation*}
h(s)=\sum_{r=0}^{n-1}b_{r}\sum_{k=0}^{\infty}\omega^{k}\sum_{n=0}^{\infty} (i\zeta)^{n} ([c;x]k)_{n}\frac{1}{n!}L^{-1}\{s^{-r-(a+b)k-\rho n-1}\}
\end{equation*}
\begin{equation}\label{2.23}
+\delta \sum_{k=0}^{\infty}\omega^{k}\sum_{n=0}^{\infty}(i\zeta)^{n}([c;x]k)_{n}\frac{1}{n!}L^{-1}\{G(s)s^{-a-(a+b)k-\rho n}\}
\end{equation}
Using $L^{-1}\{\frac{1}{s^{n}}\}=\frac{t^{n-1}}{(n-1)!}$ in the above expression, then after some calculation, we arrive at
\begin{equation*}
=\sum_{r=0}^{n-1} b_{r}\sum_{k=0}^{\infty}\omega^{k}\mu^{r+(a+b)k}E^{[c; x]k}_{\rho,1+(a+b)k+r}(i\zeta \mu^{\rho})
\end{equation*}
\begin{equation}\label{2.24}
+\delta\int_{0}^{\mu}\sum_{k=0}^{\infty}\omega^{k} (\mu-t)^{a+(a+b)k-1}E^{[c;x]k}_{\rho, a+(a+b)k}(i\zeta (\mu-t)^{\rho})g(t)dt.
\end{equation}
By using (\ref{1.15}) in (\ref{2.24}), we can achieve (\ref{2.19}) and (\ref{2.21}).
\end{proof}
\smallskip

{\bf{Funding.}} There is no source of funding for this article.\\
{\bf{Ethical approval.}} This article does not contain any studies with human participants or animals.\\
{\bf {Availability of data and materials}} Not applicable.\\
{\bf{Conflict of interest.}} The author declares that this paper involves no conflict of interest.
 
\end{document}